\documentclass[11pt]{amsart}
\usepackage{cite}

\usepackage{amsmath}
\usepackage{amssymb}
\usepackage{newlfont}
\usepackage{amsbsy}
\usepackage{mathrsfs}
\usepackage{color}
\usepackage{graphics}
\usepackage{amsthm}
\usepackage{amsfonts}
\usepackage{enumerate}

\theoremstyle{plain}
\newtheorem{lemma}{Lemma}[section]
\newtheorem*{theorem*}{Theorem}
\newtheorem*{lemma*}{Lemma}
\newtheorem*{proposition*}{Proposition}
\newtheorem*{conjecture*}{Conjecture}
\newtheorem*{corollary*}{Corollary}
\newtheorem*{problem*}{Problem}
\newtheorem{theorem}[lemma]{Theorem}
\newtheorem{conjecture}[lemma]{Conjecture}
\theoremstyle{definition}

\newtheorem{example}[lemma]{Example}
\newtheorem{remark}[lemma]{Remark}
\theoremstyle{plain}
\newtheorem{corollary}[lemma]{Corollary}
\newtheorem{proposition}[lemma]{Proposition}

\newcommand{\F}[1]{\mathscr{#1}} 

\newcommand{\R}{\mathbb{R}}

\newcommand{\Z}{\mathbb{Z}}

\newcommand{\C}{\mathbb{C}}

\newcommand{\OO}{\mathcal{O}}

\newcommand{\I}{\F{I}}

\DeclareMathOperator{\pdim}{\mathrm{pdim}}

\newcommand{\Pro}{\mathbb{P}}

\renewcommand{\P}{\Pro}

\DeclareMathOperator{\edim}{edim}

\DeclareMathOperator{\te}{\otimes}

\DeclareMathOperator{\vdim}{vdim}

\DeclareMathOperator{\Aut}{Aut}

\DeclareMathOperator{\Pic}{Pic}

\DeclareMathOperator{\Spec}{Spec}

\newcommand{\fto}[1]{\stackrel{#1}{\to}}

\newcommand{\sm}{\setminus}

\author{Jack Huizenga}

\usepackage{amsbsy}

\begin{document}

\date{\today}
\address{Department of Mathematics\\Harvard University, Cambridge, MA 02143}
\email{huizenga@math.harvard.edu}
\thanks{This material is based upon work supported under a National Science Foundation Graduate Research Fellowship}

\title{Interpolation on surfaces in $\P^3$}

\begin{abstract}
Suppose $S$ is a surface in $\P^3$, and $p_1,\ldots,p_r$ are general points on $S$.  What is the dimension of the space of sections of $\OO_S(e)$ having singularities of multiplicity $m_i$ at $p_i$ for all $i$?  We formulate two natural conjectures which would answer this question, and we show they are equivalent.  We then prove these conjectures in case all multiplicities are at most $4$.
\end{abstract}

\maketitle

\newcommand{\spacing}[1]{\renewcommand{\baselinestretch}{#1}\large\normalsize}
 \setcounter{tocdepth}{2} \tableofcontents

\section{Introduction}

Let $S\subset \P_\C^3$ be a surface of degree $d$, and $p$ a point on $S$.  The \emph{fat point} $p^m$ of multiplicity $m$ supported at $p$ is the scheme defined by the $m$th power of the ideal of the point.  Now let $p_1,\ldots,p_r$ be a general collection of points on $S$.  Given multiplicities $m_1,\ldots,m_r$, we ask the question, when does the fat point scheme $$\Gamma = p_1^{m_1}\cup \cdots\cup p_r^{m_r}$$ impose the expected number of conditions on global sections of $\OO_S(e)$?

To fix notation, we write $\F L_e^S(\Gamma)$ or  $\F L_e^S(m_1,\ldots,m_r)$ for the linear series of members of $H^0(\OO_S(e))$ containing $\Gamma$.  When some multiplicities are repeated, we may use exponential notation, so that for instance $\F L_e^S(3^a,2^b)$ denotes the series of curves having $a$ triple points and $b$ double points.  By the \emph{dimension} of a linear series we always mean the vector space dimension.  We define the \emph{virtual dimension} of such a series to be the quantity $$\vdim \F L_e^S(m_1,\ldots,m_r) = h^0(\OO_S(e)) - \deg \Gamma = h^0(\OO_S(e)) - \sum_{i=1}^r {m_i +1 \choose 2},$$ and we define the \emph{expected dimension} by $$\edim \F L_e^S(m_1,\ldots,m_r)=\max\{\vdim \F L_e^S(m_1,\ldots,m_r),0\}.$$ We call $\F L_e^S(\Gamma)$ \emph{nonspecial} if its dimension equals the expected dimension, and \emph{special} otherwise.

The case where $d=1$ so that $S=\P^2$ has received an enormous amount of attention.  The Segre-Harbourne-Gimigliano-Hirschowitz (SHGH) conjecture (see \cite{Segre}, \cite{Harbourne}, \cite{gimi2}, \cite{Hirschowitz}, respectively) states that if $\F L_e^{\P^2}(\Gamma)$ is special, then there is a multiple $(-1)$-curve $C$ in the base locus of the series.  By this we mean that the proper transform of $C$ on the blowup of $\P^2$ at the points $p_1,\ldots,p_r$ is a smooth rational curve of self-intersection $-1$.  This conjecture gives a simple algorithm for determining whether a given series $\F L_e^{\P^2}(\Gamma)$ is special or not.  For a nice history of this problem and some of the approaches used to attack it, we refer the reader to the survey article \cite{BocciMiranda} of Bocci and Miranda.  In particular, we remark that the conjecture is known to hold so long as all the multiplicities $m_i$ are bounded by at most $7$ (see Yang \cite{Yang}) or $11$, as shown in Dumnicki \cite{DumnickiJarnicki}.

Many authors have considered the problem of generalizing SHGH type statements to surfaces other than $\P^2$.  In this direction, the most common objects of study have been the rational ruled Hirzebruch surfaces $\mathbb F_n.$  Recent work of Laface \cite{Laface3} and Dumnicki \cite{Dumnicki} shows that a natural analog of the conjecture can be formulated on Hirzebruch surfaces, and that this conjecture is true when the multiplicities are small.

The SHGH conjecture also applies in the cases $d=2,3$, at least if $S$ is general.  If $S$ is a smooth quadric surface, then the blowup of $S$ at one point is isomorphic to the blowup of $\P^2$ at two points, so sections of $\OO_S(e)$ vanishing along $\Gamma$ are in bijective correspondence with sections of $\OO_{\P^2}(2e)$ vanishing along the fat point scheme $\Gamma\cup q_1^e\cup q_2^e$, and the conjecture addresses the latter case.  Likewise, thinking of a general smooth cubic surface $S$ as the blowup of $\P^2$ at six general points, sections of $\OO_S(e)$ vanishing along $\Gamma$ correspond to sections of $\OO_{\P^2}(3e)$ vanishing along $\Gamma \cup q_1^e \cup\cdots \cup q_6^e.$

Thus the first case which has not received much attention is the case $d=4$.
In this case $S$ is a $K3$-surface.  If we additionally assume $S$ is very general, so that $\Pic S = \Z$, then it has been conjectured by De Volder and Laface \cite{Laface} that the only special fat point linear series $\F L_e^S(\Gamma)$ on $S$ are the series $\F L_e^S(2e)$ for $e\geq 2$.  In particular, they show this conjecture would follow from the statement that all special linear series on $S$ have a nonreduced curve in their base locus.  It is easy to see that these series are in fact special, since the expected dimension is $0$, but taking $e$ copies of a tangent plane to $S$ produces a curve in $\F L_e^S(2e)$.

For degree $d\geq 5$, less is known.  For the same reason as in the case $d=4$, the series $\F L_2^S(4)$, $\F L_3^S(6),$ and $\F L_4^S(8)$ are special for every $d$, while $\F L_5^S(10)$ is special for $d=5$.  It seems reasonable to suspect that these are the only special series on a very general surface $S$.  In fact, we show that if every special linear series on $S$ has a multiple curve in its base locus, then these are the only special series.

\newtheorem*{thmA}{Theorem A}

\begin{thmA}
Let $S\subset \P^3$ be a very general surface of degree $d\geq 4$.  The following two statements are equivalent:
\begin{enumerate}
\item If the general member of a series $\F L_e^S(\Gamma)$ of fat points is reduced, then the series is nonspecial.

\item The only special series $\F L_e^S(\Gamma)$ are the series \begin{itemize} \item $\F L_e^S(2e)$ for $d=4$ and $e\geq 2$,
\item $\F L_e^S(2e)$ for $d=5$ and $2\leq e \leq 5$, and
\item $\F L_e^S(2e)$ for $d\geq 6$ and $2\leq e \leq 4$.
\end{itemize}
\end{enumerate}
\end{thmA}

By a very general surface, we mean a surface lying outside a countable union of closed subvarieties of the projective space of all surfaces of degree $d$.  The primary reason for assuming $S$ is very general is that then the Noether-Lefschetz theorem \cite{GriffithsHarris} implies $\Pic S=\Z$, with generator $\OO_S(1)$.  In fact, this is the only hypothesis we need to ensure that the first statement implies the second.  We refine this statement and give most of its proof in Section \ref{TwoConjsSection}.

The primary goal of the rest of this paper is to gather evidence for either of the two equivalent statements in the preceding theorem.  To do this, we study the cases where $d\geq 4$ and the multiplicities of the points are relatively small.  As in Theorem A, we concentrate on what happens when $S$ is general; it is possible to give some results when $S$ is a specific surface, but the techniques involved are much messier.
We denote by $\F L^d_e(\Gamma)$ the series $\F L_e^S(\Gamma),$ where $S$ is a general surface of degree $d$.  

\newtheorem*{thmB}{Theorem B}

\begin{thmB}
If $d\geq 4$, then the only special linear series $\F L_e^d(4^a,3^b,2^c)$ is $\F L_2^d(4)$.
\end{thmB}

It is worth remarking that if $S$ is very general of degree at least $4$, then the only line bundles on $S$ are the bundles $\OO_S(e)$.  Thus if $S$ is very general and $L$ is any line bundle on $S$, the only way fat points of multiplicity at most $4$ can fail to impose the expected number of conditions on sections of $L$ is if we are looking at the series $\F L_2^d(4)$.

We prove Theorem B by allowing a surface of degree $d$ to degenerate in a pencil to a union of two surfaces of smaller degree.  The total space of this family has singularities which we must resolve.  After resolving the singularities, we are able to modify the line bundle $\OO(e)$ on the special fiber.  We specialize some of our fat points onto each surface, and argue by induction on the degree of the surface that the series on the special fiber is nonspecial.  Then by semicontinuity, the series on the general fiber is nonspecial.  This general strategy is reminiscent of the degeneration techniques developed by Ciliberto and Miranda to study fat point series on $\P^2$ in \cite{CM1} and \cite{CM2}.  We study this degeneration in Section \ref{degenSection}.  As a bonus, this degeneration technique will allow us to complete the proof of Theorem A.

This induction on the degree of the surface will naturally lead to surfaces of degree smaller than $4$.  In this case there are many well-known techniques for addressing the question, although we will have to rely on computers to check the enormous number of cases that will arise.  After discussing what happens in these cases in Section \ref{HHSection}, we will prove Theorem B in Section \ref{QuadrupleSection}.

I would like to thank Joe Harris for the many extremely helpful discussions that led to this paper.  I would also like to thank Marcin Dumnicki for his careful critical reading of this paper, for vastly improving the argument in Section \ref{HHSection}, and for offering the proof and performing the computer calculations necessary to verify Theorem \ref{DumnickiThm}.  Finally, Ciro Ciliberto, Brian Harbourne, Antonio Laface, and Rick Miranda have also provided useful comments and guidance on this work.

\section{Two conjectures on special linear series}\label{TwoConjsSection}

The goal of this section is to establish the equivalence of two conjectures relating to special linear series of fat points on a very general surface of degree $d\geq 4$.  The first conjecture can be seen as an analog of Segre's conjecture for $\P^2$.

\begin{conjecture}\label{SegreConjS}
Let $S\subset \P^3$ be a very general surface of degree at least $4$.  If  the general member of $\F L = \F L_e^S(m_1,\ldots,m_r)$ is reduced, then $\F L$ is nonspecial.
\end{conjecture}

On the other hand, we can make another conjecture by listing all the linear series we suspect are special.  This conjecture can be viewed as an equivalent of the Harbourne-Gimigliano-Hirschowitz conjecture for $\P^2$.

\begin{conjecture}\label{HHConjS}
Let $S\subset \P^3$ be a very general surface of degree $d$ at least $4$, and let $\F L = \F L_e^S(m_1,\ldots,m_r)$ be a general series of fat points.

\begin{enumerate}
\item The series $\F L$ is special if and only if it is of the form
\begin{itemize}
\item $\F L_e^S(2e)$ if $d=4$ and $e\geq 2$,
\item $\F L_e^S(2e)$ if $d=5$ and $2\leq e\leq 5$, or
\item $\F L_e^S(2e)$ if $d\geq 6$ and $2\leq e\leq 4$.
\end{itemize}

\item If $\F L$ is nonspecial and nonreduced, then it either equals $\F L_2^S(2^3)$ (with $d$ arbitrary) or $\F L_5^S(10)$ with $d\geq 6$.

\item In all other cases, if $\F L$ is nonempty then the general member of $\F L$ is reduced and irreducible.

\item If $\F L$ is nonempty the general member of $\F L$ has exactly the assigned multiplicities at the specified base points.
\end{enumerate}
\end{conjecture}

The next result clearly implies Theorem A from the introduction.

\begin{theorem}\label{equivConjs} Conjecture \ref{SegreConjS} and Conjecture \ref{HHConjS} are equivalent.
\end{theorem}

This theorem is essentially an analog of the equivalence of B. Segre's conjecture and the Harbourne-Gimigliano-Hirschowitz conjecture for $\P^2$, shown by Ciliberto and Miranda in \cite{CM3}.  Our proof is reminiscent of the argument given there. 

We will first occupy ourselves with the more difficult direction of the theorem, showing that Conjecture \ref{SegreConjS} implies Conjecture \ref{HHConjS}.  We note that the case $d=4$ is handled in \cite{Laface}, since then $S$ is a $K3$ surface.  We thus restrict ourselves to the cases $d\geq 5$.  (However, we note that the opposite implication in case $d=4$ is not quite as straightforward as De Volder and Laface suggest in \cite{Laface}.  They later proved this implication in \cite{Laface2}.)

Denote by $S'$ the blowup of $S$ at $r$ general points $p_1,\ldots,p_r$.  Since $S$ is very general, $\Pic S$ is generated by the class $H$ of a hyperplane, and $\Pic S'$ is generated by the pullback $H$ of the hyperplane class together with the classes $E_1,\ldots, E_r$ of the exceptional divisors.

Suppose we are given an $\R$-divisor $aH - b_1E_1-\cdots - b_rE_r$ on $S'$ with $a$ a positive integer and $b_i\geq 0$ positive reals.  We define $$v(aH-\sum b_iE_i) = h^0(\OO_S(a))-1 - \sum \frac{b_i(b_i+1)}2.$$ We caution that this definition is \emph{not} the same as the definition of the virtual dimension in the sense of the rest of the paper--for the other sections of the paper, we are concerned primarily with empty linear series, for which it is convenient to use the vector space dimension, while in this section we care primarily about individual fixed curves, for which it is convenient to use the projective dimension.  To avoid any serious conflict, we will never call $v$ the virtual dimension.  When we wish to talk about the projective dimension of a series $\F L$, we will write $\pdim \F L$.

Given an effective divisor $D$ on $S$, we also write $v(D)$ for $v$ of the proper transform of $D$. So long as no confusion is likely to arise, we also denote this proper transform and its divisor class in $\Pic S'$ by $D$. Our primary tool in showing Conjecture \ref{SegreConjS} implies Conjecture \ref{HHConjS} is the following proposition (which does not depend on either conjecture).

\begin{proposition}\label{curveSumProp}
Suppose that $D$ and $D'$ are effective divisors on $S$ with $v(D)=v(D')=0$, where $d\geq 5$.  Then either \begin{enumerate}
\item $v(D+D')>0$,
\item $D=D'=H-2E_i$ for some $i$, or
\item $D=D'=H-E_i-E_j-E_k$ for some $i,j,k$.
\end{enumerate}
\end{proposition}

We will frequently apply this proposition in the form of the following immediate corollary.

\begin{corollary}\label{cont} If $D$ and $D'$ are distinct effective divisors on $S$, it is impossible to have $$v(D)=v(D')=v(D+D')=0.$$
\end{corollary}

We first single out an important lemma in the proof of the proposition.

\begin{lemma}
Consider for each positive integer $a$ the set $$X_a = \{D= aH - \sum b_iE_i:b_i\in \R_{\geq 0}, v(D)=0\}\subset N^1(S')_\R,$$ and define for each pair of positive integers $a,a'$ a function \begin{eqnarray*} X_a\times X_{a'} &\to &\R\\ (D,D')&\mapsto &v(D+D').\end{eqnarray*} This function achieves its minimum value on $X_a\times X_{a'}$ exactly at the $r$ pairs of divisors $$(D,D') = (aH-b_iE_i,a'H-b_i'E_i), \qquad (i=1,\ldots,r)$$ where $b_i$, $b_i'$ are determined by the constraints $v(D) = v(D') = 0$.
\end{lemma}

The proof of the lemma is an elementary exercise in multivariable calculus and Lagrange multipliers.  According to the lemma, to produce a lower bound for the quantity $v(D+D')$, we may examine the worst case scenario where $D$ and $D'$ are both given by a single big fat point (where the multiplicity is a real number instead of an integer).

\begin{proof}[Proof of Proposition \ref{curveSumProp}]
Define two functions $$f(a) = h^0(\OO_S(a))-1 \qquad g(a) = \frac{-1+\sqrt{1+8f(a)}}2.$$ The number $g(a)$ is the unique nonnegative number with $$\frac{g(a)(g(a)+1)}2 = f(a).$$ Say $D=aH-\sum b_i E_i$ and $D' = a'H-\sum b_i'E_i$, with $a,a'>0$.  Then by the lemma $$v(D+D')\geq v((a+a')H-(g(a)+g(a'))E_1) = f(a+a')-\frac{(g(a)+g(a'))(g(a)+g(a')+1)}2.$$ Since $$f(a+a') = \frac{g(a+a')(g(a+a')+1)}2,$$ if we are able to show that \begin{equation}\label{superadditive} g(a)+g(a')<g(a+a')\end{equation} then it will follow that $v(D+D')>0$.

Inequality \ref{superadditive} doesn't quite always hold (indeed, it can't hold in the exceptional cases of the proposition), but it holds often enough to be very useful.  To determine for which pairs $a,a'$ it does hold, assume $a\geq a'$ and write $$g(a') = (g(a')-g(a'-1))+(g(a'-1)-g(a'-2))+\cdots + (g(1)-g(0)),$$ noting that $g(0)=0$.  Suppose we have shown \begin{equation}\label{consecutiveDif}g(k+1)-g(k)>g(k)-g(k-1)\end{equation} for all integers $k\geq 2$, and that additionally $g(a+1)-g(a) > g(1)-g(0) = 2$.  It then follows that $$g(a') < (g(a+a')-g(a+a'-1))+(g(a+a'-1)-g(a+a'-2))+\cdots + (g(a+1)-g(a)),$$ from which Inequality \ref{superadditive} would follow immediately.  We note that since $d\geq 5$, we can calculate $$g(4)-g(3) \approx 2.08>2,$$ so (assuming Inequality \ref{consecutiveDif} holds) the inequality $g(a+1)-g(a)>2$ is satisfied as soon as $a\geq 3$.

To prove Inequality \ref{consecutiveDif} holds for $k\geq 2$, first observe directly that it holds for $k=2$ since $$g(3)-g(2) \approx 1.91 > 1.77 \approx g(2)-g(1).$$ For all other $k$, we write $$2(g(k+1)-g(k))= \sqrt{1+8f(k+1)}-\sqrt{1+8f(k)},$$ so Inequality \ref{consecutiveDif} amounts to showing $$\sqrt{1+8f(k+1)}+\sqrt{1+8f(k-1)}> 2\sqrt{1+8f(k)}.$$ To show this inequality, it suffices to find a convex function $G:[k-1,k+1]\to \R$ with $G(x) = \sqrt{1+8f(x)}$ for $x\in \{k-1,k,k+1\}.$  

First suppose $3\leq k\leq d-2$.  Then for $x\in\{k-1,k,k+1\}$, we
have $$f(x) = {x+3\choose 3}-1 = \frac{(x+3)(x+2)(x+1)}6-1.$$ We
thus define $$G(x) = \sqrt{1+8\left({x+3\choose 3}-1\right)}.$$ It
is easy to calculate the second derivative $G''(x)$, and we observe
that it is positive so long as $x\geq 2$.  Thus $G$ is convex on the
domain $[2,d-1]$ of interest.

On the other hand, suppose $k\geq d-1$.  In this case for $x\in \{k-1,k,k+1\}$ we have $$f(x) = {x+3\choose 3}-{x-d+3\choose 3}-1.$$ Defining $G$ in the obvious way, we can again make a straightforward calculus calculation to show $G$ is convex,  so Inequality \ref{consecutiveDif} holds for $k\geq 2$.

We have now shown Inequality \ref{superadditive} holds so long as $a\geq a'$ and $a\geq 3$.  It also holds by direct calculation in case $a=a'=2$, so the only cases remaining are $(a,a')=(2,1)$ and $a=a'=1$; in these cases Inequality \ref{superadditive} does not hold.  In case $(a,a')=(2,1)$, we have  \begin{eqnarray*}D&\in&\{\F L^S_2(3,2),\F L_2^S(2^3),\F L_2^S(3,1^3), \F L_2^S(2^2,1^3), \F L_2^S(2,1^6),\F L_2^S(1^9)\}\\ D' &\in &\{\F L_1^S(2), \F L_1^S(1^3)\}.\end{eqnarray*} One quickly checks that in every case $v(D+D')\geq 1$, no matter which points $p_i$ the multiplicities are assigned to.  On the other hand, if $a=a'=1,$ we easily see that if $D\neq D'$ then $v(D+D')>0$, but if $D=D'$ then $v(2D)\leq 0$.
\end{proof}

With the proposition in hand, the proof that Conjecture \ref{SegreConjS} implies Conjecture \ref{HHConjS} is relatievly easy.  The next results follow a similar framework to the analogous results for $K3$ surfaces given in \cite{Laface}.  We begin with a previously known result that also holds independently of the conjectures.

\begin{lemma}[Ciliberto and Chiantini \cite{CC}, Proposition 2.3]\label{exceptionalBaseLocus}
If $D=aH-\sum_i b_iE_i$ is a divisor on $S'$ with $a,b_i\geq 0$ and $|D|$ is nonempty, then $E_i$ does not lie in the base locus of $|D|$.
\end{lemma}

\begin{lemma}\label{irreducible}
Assume Conjecture \ref{SegreConjS} is true.  Suppose $\F L = \F L_e^S(\Gamma)$ is a linear series of fat points with no multiple fixed components, where $d\geq 5$.  Then the general member of $\F L$ is irreducible.
\end{lemma}
\begin{proof}
The hypothesis implies $\F L$ is nonspecial.  Blow up $S'$ at $v(\F L)$ additional general points to get a surface $S''$, and look at the member $D$ of $\F L$ passing through these points.  On $S''$ we have $|D| = \{D\}$, and $v(D)=0$.  If $D = D' + D''$ is reducible, then $v(D) = v(D') = v(D'')=0,$ contradicting Corollary \ref{cont} since $D'$ and $D''$ are distinct (observe that $D$ does not contain $E_i$ by the previous lemma).  Thus the general member of $\F L$ is irreducible.
\end{proof}

\begin{proposition}\label{fixedComps}
Assume Conjecture \ref{SegreConjS} is true.  If $d\geq 5$ and $\F L = \F L_e^S(\Gamma)$ has a multiple fixed component, then $\pdim \F L=0$.  Also, the unique member of $\F L$ has irreducible support.
\end{proposition}
\begin{proof}
Let $D$ be a member of $\F L$, and write $$|D| = \sum_{i=1}^a \mu_i C_i +\sum_{i=1}^b F_i + |D'|,$$ where the $C_i$, $F_i$ are all distinct, reduced and irreducible, $\mu_i\geq 2$ for each $i$, and $|D'|$ has no fixed components.  By Lemma \ref{exceptionalBaseLocus}, none of the $C_i$, $F_i$ are the exceptional divisors of $S'$.  Suppose $a+b\geq 2$, so that there are at least $2$ distinct fixed components; call them $C$ and $C'$.  Since $C$, $C'$, and $C+C'$ are all nonspecial and don't move in the series $|D|$, we have $v(C)=v(C')=v(C+C')=0$, contradicting Corollary \ref{cont}. Since $\F L$ has a multiple fixed component, we conclude $a=1$ and $b=0$, so $$|D| = \mu C + |D'|.$$ This also implies that if $\pdim \F L=0$, then the unique member of $\F L$ has irreducible support.

If $\pdim \F L>0$, then $\pdim |D'|>0$ and $D'$ is nonspecial by the conjecture.   Blow up $S'$ at $v(D')$ general points to get a surface $S''$; there is a unique member $F$ of $|D'|$ passing through these points.  We have $v(F)=v(C)=v(F+C)=0$ on $S''$, and $F$ is distinct from $C$, contradicting the corollary.  Thus $\pdim \F L = 0$.
\end{proof}

\begin{theorem}
Conjecture \ref{SegreConjS} implies Conjecture \ref{HHConjS}.
\end{theorem}
\begin{proof}
If $\F L$ has a multiple fixed component $C$ with multiplicity $\mu\geq 2$, then $\F L$ consists solely of $\mu C$ by Proposition \ref{fixedComps}.  Since $\mu C$ is fixed in $\F L$, we find $v(C)=0$ and $v(2C)\leq 0$, so by Proposition \ref{curveSumProp} $C$ is the sole member of either $\F L_1^S(2)$ or $\F L_1^S(1^3)$. Part (1) of Conjecture \ref{HHConjS} then follows by determining when $v(\mu C)<0$ in each case;  part (2) follows by determining when $v(\mu C)= 0$.  Part (3) is Lemma \ref{irreducible}, and part (4) is Lemma \ref{exceptionalBaseLocus}.
\end{proof}

To complete the proof of Theorem \ref{equivConjs}, we need to show that each of the special series listed in Conjecture \ref{HHConjS} in fact have a multiple curve in the base locus.  In case $d=4$, we recall the following result.

\begin{proposition}[De Volder and Laface \cite{Laface2}]\label{K3SpecialSeries}
Suppose $S$ is very general of degree $4$, and consider the series $\F L_e^S(m)$.  If $m=2e$ with $e\geq 2$, then this series is special, and its only member is the curve $eC$, where $C$ is the unique member of $\F L_1^S(2)$.  For all other $m$, this series is nonspecial.
\end{proposition}

The remaining cases will require some specialization techniques that we will develop over the course of the rest of the paper.  We will complete the proof in Section \ref{degenSection}.

\section{Computing limit linear series}\label{LimitSection}

In this section, we discuss an important tool for proving the nonspeciality of linear series.  We refer the reader to Evain \cite{Evain} and Roe \cite{Roe} for more formal and general expositions of the technique, and also for the proof of the main theorem in this section.  

Let $S$ be a surface with a line bundle $L$, and let $V\subset H^0(L)$ be a linear series.  Let $C\subset S$ be an irreducible curve, and let $p_t^m$ be a family of fat points of multiplicity $m$, parameterized by a disk $\Delta$, approaching a general point $p_0$ of $C$ transversely as $t\to 0$.

If $\Gamma\subset S$ is a zero-dimensional scheme, we denote by $V(-\Gamma)$ the subseries of members of $V$ containing $\Gamma$.  On the other hand, if $D\subset S$ is an effective divisor, we denote by $V(-D)$ the subseries of divisors $E$ in $H^0(L(-D))$ such that $E+D\in V$.

The dimension $\dim V(-p_t^m)$ is an upper-semicontinuous function of $t$.  Thus for generic $t$ it assumes some minimal value $v$.  A priori (and in many cases of interest) this dimension jumps for $t=0$.  However, since the Grassmannian $G(v,V)$ is proper, we can ask what the limit $$V_0 = \lim_{t\to 0} V(-p_t^m)$$ is.  At any rate, it will be a $v$-dimensional subspace of $V(-p_0^m)$.  However, if we can find explicit geometric conditions describing $V_0$, it may be possible to compute the dimension of $V_0$, and hence compute the dimension of $V(-p_t^m)$ for general $t$.

\begin{example} For a basic example, let us look at $S=\P^2$.  Take $C=L$ to be a line, and let $V\subset H^0(\OO_{\P^2}(2))$ be the $4$-dimensional series of conics passing through two points $q,r$ lying on $L$.  Now let $p_t^2$ be a double point tending to $L$ as $t\to 0$.  Clearly $V(-p_0^2)$ is $2$-dimensional, consisting of all reducible conics containing $L$ and singular at $p_0$.  We have ``lost'' a condition, since the line $L$ meets $p_0^2\cup q \cup r$ in a subscheme of degree $4$, while it only takes $3$ conditions to force $L$ to appear in the base locus: $h^0(\OO_L(2)) = 3$.

The solution to this problem is to observe that we want $p_t^2$ to meet $L$ in a subscheme of degree $1$, not of degree $2$.  Of course, this is impossible:  either $t\neq 0$, in which case $p_t^2$ and $L$ are disjoint, or $t=0$, in which case $p_t^2$ meets $L$ in a subscheme of degree $2$.  However, if we  restrict our family to the nonreduced base $\Spec k[t]/(t^2)$ then we find that ``for general $t$'' $p_t^2$ will contain the subscheme $p_0$ of $L$, without containing the subscheme $2p_0\subset L$.

Thus when $p_t^2$ ``hovers'' in a first order neighborhood of $L$, we find that $p_t^2 \cup q\cup r$ meets $L$ in a scheme of degree $3$.  It follows that $L$ lies in the base locus of $V(-p_t^2)$.  To remove it from the series, we have to look at members of $V(-L)$ containing the scheme $(p_t^2:L)$ defined by the ideal quotient of $\I_{p_t^2}$ by $\I_L$ (where the ideal quotient is calculated keeping $t^2=0$ in mind!).  This ideal quotient is a bit nasty, but upon setting $t=0$ (i.e. restricting from our nonreduced base to $\Spec k[t]/(t)$), it just defines the subscheme $2p_0\subset L$.  In other words, members of $V_0(-L)$ are all tangent to $L$ at the point $p$, which is to say that $V_0$ consists solely of the double line $2L$.
\end{example}

For the purposes of the present paper, we will only need a mild generalization of the previous example.  In the notation of the second paragraph of this section, let $v$ be the dimension of $V|_C$, and assume $v\leq m$, so that we expect $C$ to appear in the base locus after specializing $p^m_t$ onto $C$.  If we can restrict $t$ in such a way that $p_t^m$ meets $C$ in a scheme of degree $v$, then $C$ will appear in the base locus of $V(-p_t^m)$ (since $p_0$, being general, is not an inflectionary point for $V|_C$), and no conditions will be ``wasted'' when we remove $C$ from the base locus.  In order to achieve this, we can restrict our family to $\Spec k[t]/(t^{m-v+1}).$  We then calculate the scheme $\Gamma$ whose ideal sheaf is given by the ideal quotient of $\I_{p_t^m}$ by $\I_C$; members of $V(-p_t^m)$ residual to $C$ must contain this scheme.  Setting $t=0$, if $(f,x)$ is a system of parameters for the maximal ideal of the local ring $\OO_{S,p_0}$ at $p_0$ and $f=0$ is a local defining equation for $C$ at $p_0$, then $$\I_{\Gamma_0} = (x^{m},x^{m-1}f,\ldots,x^{v+1}f^{m-v-1} ,x^{v-1}f^{m-v},x^{v-2}f^{m-v+1},\ldots,f^{m-1}).$$
Intuitively, every member of $V_0(-C)$ must have a singularity of mutliplicity at least $m-1$ at $p_0$, and the tangent cone of the singularity must contain the tangent line to $C$ at $p_0$ with multiplicity $m-v$ (unless the curve actually has a singularity of multiplicity $m$).  Thus the $m-v$ extra conditions which would have been ``lost'' by specializing $p_t^m$ naively appear as extra tangency conditions on the branches of the singularity of the limit curves.

It is worth noting that the choice of curve $\Delta\to S$ along which $p_t^m$ approaches $p_0\in C$ is essentially irrelevant to the limit computation, so we will usually not discuss it.

\begin{theorem}[Evain \cite{Evain}]\label{LimitTheorem}
Let $S$ be a surface, let $L$ be a line bundle on $S$, and let $V\subset H^0(L)$ be a linear series.  Let $C\subset S$ be an irreducible curve, and put $v=\dim V|_C$.   Let $m\geq v$, and let $p^m$ be a general fat point of multiplicity $m$ on $S$.  Let $$V_0 = \lim_{p\to p_0} V(-p^m).$$ Then $V_0$ contains $C$ in its base locus, and every member of $V_0(-C)$ either 
\begin{enumerate}\item has a singularity at $p_0$ of multiplicity $m-1$, with tangent cone containing the tangent line to $C$ with multiplicity $m-v$, or
\item has a singularity of multiplicity at least $m$.
\end{enumerate}
\end{theorem}

To avoid having to repeat the geometric conclusion of the theorem, we define schemes 
$$\delta_{m,n} = \Spec\C[x,y]/(x^{m+1},x^my,\ldots,x^{m-n+2}y^{n-1},x^{m-n}y^n,x^{m-n-1}y^{n+1},\ldots,y^m) \qquad (n\leq m).$$
We refer to a closed subscheme of $S$ isomorphic to $\delta_{m,n}$ as a $\delta_{m,n}$-point.   If $n>0$, then at any smooth point $p$ of $S$ there is a unique $\delta_{m,n}$-point supported at $p$ and ``pointing'' in a given tangent direction.  A curve containing a $\delta_{m,n}$-point has a singularity of multiplicity $m$ and has tangent cone containing the distinguished direction with multiplicity $n$ (unless it has a singularity of multiplicity $m+1$).  In this notation, the theorem says  $$V_0(-C) \subset V(-C)(-\delta_{m-1,m-v})$$ 

Since it is easier to deal with ordinary fat points than the schemes $\delta_{m,n}$ (for one thing, there is a unique fat point supported at a given point $p\in S$, as opposed to the $1$-parameter family of $\delta_{m,n}$-points), it is useful to be able to reduce questions involving $\delta_{m,n}$-points to questions only involving fat points.  Given a series $V$, we would expect that a general $\delta_{m,n}$-point on $S$ imposes $\min\{\deg \delta_{m,n},\dim V\}$ conditions on curves in $V$; we note that $$\deg \delta_{m,n} = {m+1\choose 2}+n.$$   The following lemma will allow us to relate many series involving $\delta_{m,n}$-points to simpler cases involving only fat points.

\begin{lemma}\label{DeltaLemma}
Let $S$ be an irreducible surface with a line bundle $L$, and let $V\subset H^0(L)$ be a linear series.  Assume a general $m$-uple point imposes the expected number of conditions on $V$.  Then a general $\delta_{m,n}$-point either imposes independent conditions on $V$ or imposes the same number of conditions as a general $(m+1)$-uple point.
\end{lemma}

\begin{proof}
Let $p$ be a smooth point of $S$ such that $V(-p^m)$ and $V(-p^{m+1})$ both have as small a dimension as possible.  Blow up $S$ at $p$, and let $E\cong \P^1$ be the exceptional divisor of the blowup.  Restricting $V(-mE)$ to $E$ yields a series $W\subset H^0(\OO_E(m))$, and we have an exact sequence $$0\to V(-(m+1)E)\to V(-mE)\to W\to 0.$$  Observe that $\dim W + {m+1 \choose 2}$ is equal to the number of conditions imposed by $p^{m+1}$ on $V$.  Now if $q$ is a general point on $E$, it is not an inflectionary point for the series $W$, and hence $W(-nq)$ has codimension $\min\{n,\dim W\}$ in $W$.  If we place a $\delta_{m,n}$-point at $p$ and pointing in the direction of $q$, then $V(-\delta_{m,n})$ is identified with the preimage of $W(-nq)$ in $V(-mE)$, and thus has codimension ${m+1\choose 2} + \min\{n,\dim W\}$ in $V$.  In case $n\leq \dim W$, we see $\delta_{m,n}$ imposes independent conditions on $V$; otherwise, it imposes the same number of conditions as $p^{m+1}$.
\end{proof}

\section{A degeneration method}\label{degenSection}

Our approach for proving the nonspeciality of fat point series on a general surface of degree $d$ will be to degenerate a surface of degree $d$ into a reducible surface; the particular degeneration we consider was studied in \cite{GriffithsHarris} to give a proof of the Noether-Lefschetz theorem without using Hodge theory.  It is then also possible to modify the line bundle on the special surface, which provides useful extra freedom.  This approach is reminiscent of the techniques used by Ciliberto and Miranda \cite{CM1} to prove the nonspeciality of series on $\P^2$; however, there are some subtle issues to deal with which don't appear when working with the plane.

To start with, fix a decomposition $d = s+t$.  Let $S$ and $T$ be smooth surfaces of degrees $s$ and $t$, meeting transversely in a smooth irreducible curve $C$, and let $U$ be a general smooth surface of degree $d$, chosen generically with respect to $S$ and $T$.  Let $X\subset \P^3\times \Delta$ be the total space of the pencil of surfaces of degree $d$ spanned by $S\cup T$ and $U$, with $X_0 = S\cup T$.

The threefold $X$ has one major defect for our purposes:  it fails to be smooth at the $dst$ points $p_1,\ldots,p_{dst}$ of the intersection $U\cap C$.  However, $U$ being general, these singularities are all ordinary double points.  When we blow them up in $X$, the exceptional divisor lying over $p_i$ is isomorphic to a nonsingular quadric surface $Q_i$ for each $i$.  This also blows up both $S$ and $T$, and their exceptional divisors form lines of opposite rulings in $Q_i$.

It is now possible to blow down the $Q_i$ along either ruling without creating new singularities in the threefold.  We blow down each $Q_i$ along the ruling containing the exceptional divisor of the blowup of $S$ at $p_i$, and call the resulting smooth threefold $\tilde X$.  The central fiber $\tilde X_0$ of $\tilde X$ is isomorphic to the union of the blowup $\tilde T$ of $T$ along the points $p_1,\ldots,p_{dst}$ with a surface isomorphic to $S$ (which we'll continue to denote by $S$);  they meet along a curve $\tilde C$ isomorphic to $C$, sitting in $\tilde T$ as the proper transform of $C$, and in $S$ as $C$ does in $S$.  Denote by $E_1,\ldots, E_{dst}$ the exceptional curves lying over $p_1,\ldots,p_{dst}$, and by $H$ the hyperplane class in $\tilde T$.

Our threefold $\tilde X$ comes equipped with a natural map $\alpha : \tilde X\to \P^3$ contracting the exceptional curves $E_i$.  The pullback $L = \alpha^\ast \OO_{\P^3}(e)$ restricts to $\OO_{X_t} (e)$ on each fiber $\tilde X_t = X_t$ for $t\neq 0$.  Similarly, it restricts to $\OO_S(e)$ on $S\subset \tilde X_0$, and to $\OO_{\tilde T}(eH)$ on $\tilde T$.

We can now modify the line bundle $L$ on $\tilde X$ by twisting by a multiple $\mu$ of the line bundle $\OO_{\tilde X}(S)$.  Clearly $\OO_{\tilde X}(S)$ is trivial on $\tilde X_t$ for time $t\neq 0$.  To analyze the situation on the special fiber, we observe that since $S\cap \tilde T = \tilde C$, we have $$\OO_{\tilde X}(S)|_{\tilde T} = \OO_{\tilde T}(\tilde C) = \OO_{\tilde T}(sH-E_1-\cdots - E_{dst}) \qquad \textrm{and}\qquad \OO_{\tilde X}(S)|_S=\OO_S(-t).$$ Thinking of the first equality in terms of fat points on the original surface $T$, we see that the global sections of $ (L\te \OO_{\tilde X}(\mu S))|_{\tilde T}$ correspond to sections of $\OO_T(d+ s\mu)$ containing fat points of multiplicity $\mu$ at each of the points $p_1,\ldots,p_{dst}$.  It is of course important to note that these points are not in general position, either on $C$ or on $T$: in fact, these $dst$ points are of the divisor class $dH \in \Pic C$.

At this point, we are ready to start incorporating fat points into the picture.  Let's let $\Gamma,\Gamma'\subset \tilde X$ each be schemes flat over $\Delta$, of relative dimension $0$, such that $\Gamma_0$ is supported in $S \sm C$ and $\Gamma'_0$ is supported in $\tilde T\sm \tilde C$.  Let $M$ be any line bundle on $\tilde X$, and suppose we are trying to show that $h^0(X_t,(M \te\I_{\Gamma\cup \Gamma'})_t)=0$ for general $t$.  By semicontinuity, $$h^0(X_t,(M\te \I_{\Gamma\cup \Gamma'})_t) \leq h^0(\tilde X_0,(M\te \I_{\Gamma\cup \Gamma'})_0)$$ for $t\neq 0$.  In an attempt to show the latter number is zero, we can write the global sections on the special fiber as a fiber product $$H^0(\tilde X_0,(M\te \I_{\Gamma\cup \Gamma'})_0) \cong H^0(S,(M\te \I_\Gamma)|_S)\times_{H^0(C,M|_C)} H^0(\tilde T,(M\te \I_{\Gamma'})|_{\tilde T}).$$ Then to show the fiber product is empty, we must show three things:
\begin{enumerate}
\item the restriction map $H^0(S,(M\te \I_{\Gamma})|_S)\to H^0(C,M|_C)$ is injective,
\item the restriction map $H^0(\tilde T,(M\te \I_{\Gamma'})|_{\tilde T})\to H^0(C,M|_C)$ is injective, and
\item the images of these two restriction maps intersect in $0$.
\end{enumerate}

In the present circumstance, we will always take $M$ to be one of the line bundles $L\te \OO_{\tilde X}(\mu S)$ for some nonnegative number $\mu$, and $\Gamma, \Gamma'$ will be moving families of general fat points, with $\Gamma$ limiting to general points of $S$ over $t=0$ and $\Gamma'$ limiting to general points of $\tilde T$ over $t=0$.  With these choices, we have $$H^0(S,(M\te \I_{\Gamma})|_S) = \F L_{e-t\mu}^S(\Gamma_0) \qquad \textrm{and}\qquad  H^0(\tilde T,(M\te \I_{\Gamma'})|_{\tilde T}) = \F L_{e+s\mu}^T(\Gamma'_0;\mu^{dst}),$$ where the $\mu$-uple points in the second series are supported on a divisor of class $dH$ on $C$; the fat points in $\Gamma_0$ and $\Gamma_0'$ are in general position (we separate $\Gamma_0'$ from $\mu^{dst}$ be a semicolon to indicate that the $dst$ points are in fact in special position).  Also, $(L\te \OO_{\tilde X}(\mu S))|_C = \OO_C(e-t\mu)$.

We will typically apply this discussion in the form of the next result.

\begin{proposition}\label{DegenProp}
Let $\Gamma$ and $\Gamma'$ be general collections of fat points on surfaces $S$ and $T$ of degrees $s,t$, respectively, with $s+t=d$.  Put $C = S\cap T$, and let $\mu$ be a nonnegative integer. Suppose that
\begin{enumerate}
\item the series $\F L^S_{e-t(\mu+1)}(\Gamma)=0,$ and
\item the subseries of members of $\F L_{e+s\mu}^T(\Gamma';\mu^{dst})$ restricting to $\F L^S_{e-t\mu}(\Gamma)|_C\subset H^0(\OO_C(e-t\mu))$ is empty.
\end{enumerate}
Then $\F L^d_e(\Gamma,\Gamma')$ is empty.
\end{proposition}
\begin{proof}
The first hypothesis ensures that property (1) from the previous discussion holds.  On the other hand, the second hypothesis ensures that properties (2) and (3) hold.
\end{proof}

To check the second condition when applying the proposition, we take the following approach.  First, find the dimension of the series $\F L_{e-t\mu}^S(\Gamma)$, typically by induction.  Since the restriction map $\F L_{e-t\mu}^S(\Gamma)\to H^0(\OO_C(e-t\mu))$ is injective by (1), this allows us to calculate the dimension of the image of the restriction map; write $W = \F L_{e-t\mu}^S(\Gamma)|_C$.  It is typically difficult to obtain any more information about $W$ than its dimension $w$.  Denote by $\F L_{e+s\mu}^T(\Gamma';\mu^{dst};w)$ the subseries of members of $\F L_{e+s\mu}(\Gamma';\mu^{dst})$ restricting to $W$.  We must show this series is empty.  To do so, we first specialize some of the points in $\Gamma'$ onto the curve $C$, so that they give at least $w$ conditions along $C$, using Theorem \ref{LimitTheorem} to calculate the residual schemes.  The important thing to note is that the residual schemes calculated by Theorem \ref{LimitTheorem} depend only on the number $w$, instead of the whole series $W$!  Thus after $C$ splits, the series we are looking at is specified entirely by the geometric data of containing certain schemes (either fat points or $\delta_{m,n}$-points, perhaps in somewhat special position).  We must then prove this residual series is empty.

\begin{remark} In previous works that made use of these kinds of degeneration arguments, it was usually the case that one of the surfaces is $\P^2$ and the curve $C$ is a line on that $\P^2$.  One then appeals to the following lemma:  given any two linear series $V,W\subset H^0(\OO_{\P^1}(d)),$ there is an automorphism $g\in \Aut \P^1$ such that $gV$ intersects $W$ transversely.  We then can realize this automorphism by an automorphism of $\P^2$ fixing the line, and moving the fat points on $\P^2$ accordingly allows us to assume the images of the restriction maps meet transversely.  

This lack of a transversality lemma is the primary reason we must use the specialization technique of Section \ref{LimitSection} instead of easier methods.
\end{remark}

To see the proposition in action, let's verify the nonspeciality of a relatively simple system that appears to be difficult to approach without the degeneration method.

\begin{lemma}\label{twoQuadruple}
The series $\F L^4_3(4^2)$ is empty.
\end{lemma}
\begin{proof}
We degenerate our surface of degree $4$ to a union $S\cup T$ of two general quadrics, meeting along an elliptic curve $C$ of degree $4$.  We take $\mu=1$ in the proposition, and specialize both quadruple points onto $T$.  The kernel system $H^0(\OO_S(-1))$ on $S$ vanishes.  On $T$, we have the series $\F L_5^2(4^2;1^{16})$, where the $16$ simple points lie along $C$ and have divisor class $4H\in \Pic C$.  Now $H^0(\OO_C(e-t\mu)) = H^0(\OO_C(1))$ is $4$-dimensional, as is the series $\F L_1^2(\emptyset)$ on $S$, so we must show that the series $\F L_5^2(4^2;1^{16})$ is empty.  Specializing a quadruple point onto $C$ causes it to split, and leaves a triple point behind.  The triple point is still in general position on $T$, since $S$ is general.  Thus we are reduced to the series $\F L_3^2(4,3)$, which is easily seen to be empty by looking at the lines in the rulings of the quadric surface.
\end{proof}

For a more substantial application of this setup, we can complete the proof that Conjecture \ref{HHConjS} implies Conjecture \ref{SegreConjS}.  The key case is taken care of in the next lemma.

\begin{lemma}\label{L55(10)}
The series $\F L_5^5(10)$ has the curve $5C$ as its sole member, where $C=H-2E$.
\end{lemma}
\begin{proof}
Degenerate a quintic surface into the union of a plane $S$ and a very general quartic surface $T$, meeting along a smooth plane quartic curve $C$.  Take $\mu=1$ in the proposition, and specialize the $10$-uple point onto $T$.  The kernel system on $S$ vanishes, and the system we're left with on $T$ is the system $\F L_6^4(10;1^{20})$.  The $20$ simple points impose only $19$ conditions on the $22$-dimensional $H^0(\OO_C(6))$; specializing the $10$-uple point onto $C$ and calculating the limit series with Theorem \ref{LimitTheorem} shows that the residual system is $\F L_5^4(\delta_{9,7})$.  Now $h^0(\OO_T(5))=52$, and we know from Proposition \ref{K3SpecialSeries} that a $10$-uple point imposes exactly $51$ conditions on $H^0(\OO_T(5))$.  Furthermore, Proposition \ref{K3SpecialSeries} shows that a $9$-uple point imposes independent conditions on $H^0(\OO_T(5))$.  By Lemma \ref{DeltaLemma} we conclude that a general $\delta_{9,7}$-point imposes $51$ conditions on $H^0(\OO_T(5))$, which is to say that the series $\F L_5^4(\delta_{9,7})$ contains a unique member.  Thus the series on the general quintic surface has a unique member, which must be $5C$.
\end{proof}

\begin{proof}[Proof of Theorem \ref{equivConjs}]
It remains to be shown that the series $\F L_e^d(2e)$ consists solely of the curve $eC$ for $2\leq e \leq 4$ and $d\geq 5$.  The case $d=5$ is clearly implied by Lemma \ref{L55(10)}.  We proceed by induction on $d$.  Degenerate a surface of degree $d$ into a plane $S$ and a general surface $T$ of degree $d-1$, and specialize the fat point onto $T$.  We take $\mu=0$, leaving the line bundle unmodified.  The series on $S$ is $H^0(\OO_S(e))$, and the restriction map $H^0(\OO_S(e))\to H^0(\OO_C(e))$ is an isomorphism since $e<d-1$.  On $T$, we have the series $\F L_e^{d-1}(2e)$, which has a unique member by induction, and this glues with a unique section on $S$ to give a unique section on the special fiber.
\end{proof}

To apply Proposition \ref{DegenProp} systematically, we'll need to know that this degeneration method preserves the virtual dimension of our series, in an appropriate sense.  By the virtual dimension of the series $\F L_{e+s\mu}^T (\Gamma';\mu^{dst};w)$ with the $dst$ points in the appropriate special position, we will mean the number \begin{eqnarray*}\vdim \F L_{e+s\mu}^T (\Gamma';\mu^{dst};w) &:=& h^0(\OO_{\tilde T}(eH + \mu \tilde C)) - \deg \Gamma'-(h^0(\OO_{C}(e-t\mu))-w)\\&=&\dim \F L_{e+s\mu}^T(\mu^{dst}) -\deg \Gamma'-(h^0(\OO_{C}(e-t\mu))-w).\end{eqnarray*}
We define the virtual dimension of a series $\F L^T_{e+s\mu}(\Gamma';\mu^{dst})$ by setting $w = h^0(\OO_C(e-t\mu))$ in the above definition. Note that this is \emph{not} the same number as if we considered $\F L_{e+s\mu}^T(\Gamma';\mu^{dst})$ as a series of fat points in general position!  This definition of the virtual dimension takes into account the fact that the points $\mu^{dst}$ are in special position, and generally do not impose independent conditions on sections of $\OO_T(e+s\mu)$.

For the particular types of degenerations that we will use in proving the quadruple point theorem, the next theorem shows that the virtual dimension of the original series is the same as the virtual dimension of the series on the surface $T$.  In particular, if the original series has nonpositive virtual dimension, so does the series on $T$.

\begin{theorem}\label{edimTheorem}
Let $\Gamma$ and $\Gamma'$ be general collections of fat points on surfaces $S$ and $T$ of degrees $s$, $t$, respectively, with $s+t = d$.  Put $C = S\cap T$, and let $\mu$ be a nonnegative integer.  Suppose
\begin{enumerate}
\item the series $\F L^S_{e-t(\mu+1)}(\Gamma)=0$, and
\item the series $\F L^S_{e-t\mu}(\Gamma)$ is nonspecial, of nonnegative virtual dimension $w$.
\end{enumerate}
Then $$\vdim \F L_e^d(\Gamma,\Gamma') = \vdim \F L_{e+s\mu}^T(\Gamma';\mu^{dst};w),$$  and $\F L_e^d(\Gamma,\Gamma')$ is empty if $\F L_{e+s\mu}^T(\Gamma';\mu^{dst};w)$ is empty.
\end{theorem}
\begin{proof}
Proposition \ref{DegenProp} implies everything except the equality of virtual dimensions.  To see the equality, let $U$ be a surface of degree $d$, and let $\tilde X\to \Delta$ be the threefold constructed earlier, corresponding to $U$ and $S\cup T$.  Unraveling the definitions, the above equality simplifies to showing $$h^0(\OO_S(e-t\mu))+h^0(\OO_{\tilde T}(e H + \mu \tilde C))-h^0(\OO_C(e-t\mu)) = h^0(\OO_U(e)).$$ Letting $M$ be the line bundle $\alpha^\ast \OO_{\P^3}(e) \te \OO_{\tilde X}(\mu S)$, we recognize the number on the left hand side as $h^0(\tilde X_0,M|_0)$ since $$H^0(\tilde X_0,M|_0)\cong H^0(\OO_S(e-t\mu))\times_{H^0(\OO_C(e-t\mu))} H^0(\OO_{\tilde T}(eH+\mu \tilde C))$$ is a fiber product and restriction $H^0(\OO_S(e-t\mu))\to H^0(\OO_C(e-t\mu))$ is surjective.  We prove the displayed equality by induction on $\mu$.  At the same time, we will show the maps $$H^0(\OO_{\tilde T}(e H+\mu \tilde C))\to H^0(\OO_C(e-t\mu))$$ are surjective by induction on $\mu$.

Write $L = M|_0$.  For $\mu = 0$, we note that $h^1(L)=0$.  Indeed, examining the exact sequence $$0\to L \to L|_S \oplus L|_{\tilde T} \to L|_C\to 0,$$ we find that $H^0(L|_S\oplus L|_{\tilde T})\to H^0(L|_C)$ is surjective since $H^0(\OO_S(e))\to H^0(\OO_C(e))$ is surjective.  Thus $H^1(L)$ injects into $H^1(L|_S)\oplus H^1(L|_{\tilde T})$, and this latter group clearly vanishes.  This means the dimension $h^0(\tilde X_t,M|_t)$ doesn't jump at $t=0$, so $h^0(L) = h^0(\OO_U(e))$.  We also note that the map $H^0(\OO_{\tilde T}(eH))\to H^0(\OO_C(e))$ is clearly surjective, since we can identify $H^0(\OO_{\tilde T}(eH))$ with $H^0(\OO_T(eH))$.

Assume our inductive hypothesis holds for $\mu-1$.  From the exact sequence $$0\to \OO_{\tilde T}(eH+(\mu-1)\tilde C)\to \OO_{\tilde T}(eH+\mu \tilde C)\to \OO_C(e-t\mu)\to 0$$ we deduce an inequality $$h^0(\OO_{\tilde T}(eH+\mu\tilde C))\leq h^0(\OO_{\tilde T}(eH+(\mu-1)\tilde C))+h^0(\OO_C(e-t\mu)),$$ with equality iff $H^0(\OO_{\tilde T}(eH+\mu \tilde C))\to H^0(\OO_C(e-t\mu))$ is surjective.  Thus \begin{eqnarray*}h^0(L)&=&h^0(\OO_S(e-t\mu)) + h^0(\OO_{\tilde T}(eH + \mu \tilde C))-h^0(\OO_C(e-t\mu))\\
&\leq & h^0(\OO_S(e-t\mu))+h^0(\OO_{\tilde T}(eH + (\mu-1)\tilde C))\\
&=& h^0(\OO_S(e-t\mu))  + h^0(\OO_C(e-t(\mu-1)))-h^0(\OO_S(e-t(\mu-1)))+ h^0(\OO_U(e))\\
&=& h^0(\OO_U(e)),
\end{eqnarray*}
with the last equality coming from the exact sequence $$0\to H^0(\OO_S(e-t\mu))\to H^0(\OO_S(e-t(\mu-1)))\to H^0(\OO_C(e-t(\mu-1)))\to 0.$$  Therefore $h^0(L)\leq h^0(\OO_U(e)),$ and by semicontinuity equality holds.
\end{proof}

\section{Double, triple, and quadruple points for $1\leq d \leq 3$}\label{HHSection}

Since the argument behind our theorem for fat points of multiplicity at most $4$ on a surface of degree at least $4$ depends on induction on the degree $d$ of the surface, we must discuss what happens when $1\leq d\leq 3$.  In all three of these cases, the SHGH conjecture predicts which series are special.

In case $d=1$ and all the multiplicities of the points are at most $4$, the SHGH conjecture is known to hold (in fact, it is known to hold so long as the multiplicities are at most $11$).  We will not actually need the case $d=1$ in our induction, so we concentrate on the cases $d=2$ and $d=3$.

In case $d=2$, a series $\F L_e^2(4^a,3^b,2^c)$ is nonspecial if and only if the series $\F L_{2e}^1(e^2,4^a,3^b,2^c)$ is nonspecial.  Similarly, for $d=3$ a series $\F L_e^3(4^a,3^b,2^c)$ is nonspecial if and only if the series $\F L_{3e}^1(e^6,4^a,3^b,2^c)$ is nonspecial.  Thus we are reduced to determining when these types of planar series are nonspecial.

A planar series $\F L_e^1(m_1,\ldots,m_r)$ with $m_1\geq m_2\geq \cdots \geq m_r$ is said to be \emph{standard} if $m_1+m_2+m_3\leq e$.   Recall that a series is called $(-1)$-\emph{special} if it has a multiple $(-1)$-curve in its base locus (and thus is also special).  The SHGH\ conjecture states that a special series in $\P^2$ is $(-1)$-special.  An important result on standard series is the following.

\begin{proposition}[Gimigliano \cite{Gimi}]
A standard series is never $(-1)$-special.  Assuming the SHGH conjecture, a standard series is nonspecial.
\end{proposition}

Now suppose we are given a series $\F L_e^1(m_1,\ldots,m_r)$.  If $m_1+m_2\geq e+1$, then there is a line in the base locus and we remove it from the system, adjusting $e,m_1,m_2$ accordingly--if this changes the expected dimension of our series, it is special.  If the series is nonstandard, put $a=m_1+m_2+m_3-e>0$.  Since $m_1+m_2\leq e$, we have $a\leq m_3$. Then by applying a Cremona transformation centered at the three points of largest multiplicity, we can transform this nonstandard series to the series  $$\F L_e^1(m_1,\ldots,m_r) \fto{\textrm{Cremona}} \F L_{e-a}^1(m_1-a,m_2-a,m_3-a,m_4,\ldots,m_r)$$ of the same dimension and expected dimension.  This decreases the degree $e$.  Continuing this process of splitting lines and performing Cremona transformations will eventually leave us with a standard system; if it has the same expected dimension as our original system then by SHGH the original system is nonspecial.  Otherwise, it is special.

\begin{proposition}\label{d23}
Assume the SHGH conjecture is true.
\begin{enumerate}
\item The only special linear series $\F L_e^2(4^a,3^b,2^c)$ are
\begin{itemize}
\item $\F L_2^2(2^3)$, $\F L_2^2(4)$,
\item $\F L_3^2(2,3^2)$, $\F L_3^2(3^3)$, $\F L_3^2(4,2^2)$,
\item $\F L_4^2(4,2^5)$, $\F L_4^2(4,3^2)$, $\F L_4^2(4,3^2,2)$, $\F L_4^2(4^2,2)$, $\F L_4^2(4^2,2^2)$, $\F L_4^2 (4^2,3)$, $\F L_4^2(4^3)$,
\item $\F L_5^2(4^3)$, $\F L_5^2(4^3,2)$, $\F L_5^2(4^3,2^2)$, $\F L_5^2(4^3,3)$, and
\item $\F L_6^2(4^5)$.
\end{itemize}
\item The only special linear series $\F L_e^3(4^a,3^b,2^c)$ is $\F L_2^3(4)$.
\end{enumerate}
\end{proposition}

We would like to thank Marcin Dumnicki for showing us the following proof.
\begin{proof}
(1) A series $\F L_{2e}^1(e^2,4^a,3^b,2^c)$ is nonstandard so long
as one of $a,b,c$ is nonzero.  If $a\neq 0$, we apply a Cremona
transformation centered at the two $e$-uple points and a quadruple
point to arrive at a series $$\F
L_{2e-4}^1((e-4)^2,4^{a-1},3^b,2^c).$$ If $e\geq 8$, then this
series is standard, hence nonspecial by the SHGH conjecture.   Now
that the problem is finite (and not very large) we can use Cremona
transformations as outlined above to determine the special series
with $a\neq 0$ and $e\leq 7$.  The cases where $a=0$ are similar.

(2) Here things are even easier.  A series $\F L_{3e}^1(e^6,4^a,3^b,2^c)$ is standard unless $e\leq 3,$ and there are then only a handful of cases to check with Cremona transformations.
\end{proof}

The above result is all we need to give a proof of Theorem C conditional on the SHGH conjecture.  However, we note that standard techniques are strong enough to prove Proposition \ref{d23} without assuming the SHGH conjecture (although the proof is far from pretty).  The proof uses the techniques of Dumnicki-Jarnicki \cite{DumnickiJarnicki} which are substantially different from those developed so far in the paper, and requires checking some 22,680 cases by computer.  We are grateful to Marcin Dumnicki for showing us the strategy and carrying out the necessary computer computation.

\begin{theorem}[Dumnicki]\label{DumnickiThm}
Proposition \ref{d23} is true without assuming the SHGH conjecture.
\end{theorem}

The basic idea of the proof is that the techniques of \cite{DumnickiJarnicki} allow one to reduce the problem to a sizable finite computation, which can then be done by computer.  We omit the proof, since it would take us far afield of the methods in this paper.

\section{Double, triple, and quadruple points on a general surface of degree $\geq 4$}\label{QuadrupleSection}

We are now ready to prove our main theorem on points of multiplicity at most $4$ on a surface of degree at least $4$.

\begin{thmB}
For $d \geq 4$, the only special systems $\F L^d_e(4^a,3^b,2^c)$ are the systems $\F L_2^d(4)$.
\end{thmB}

The proof is primarily by induction on $d$, with the case $d=4$ being far more difficult than the others.  The primary difficulty for $d=4$ is that if we degenerate to a union of two quadric surfaces, there are lots of special linear series on both quadrics (in light of Proposition \ref{d23}).  On the other hand, for $d\geq 5$, we can always degenerate the surface to a union of two surfaces where there are only lots of special series on one of the two surfaces.

\begin{proof}
Consider a series $\F L = \F L_e^d(4^a,3^b,2^c)$ different from $\F L_2^d(4)$.  By adding simple points to $\F L$, we may assume that $\vdim \F L\leq 0$, and we must show $\F L$ is empty.  While there are several cases to consider, the basic approach is the same in each case.  We must degenerate our surface of degree $d$ into a union $S\cup T$ of surfaces of degree $s$ and $t$, meeting along a curve $C$.  We modify the line bundle on $S\cup T$ by choosing the parameter $\mu$ in Theorem \ref{edimTheorem}, and we specialize some of our points $\Gamma$ onto $S$ and the rest $\Gamma'$ onto $T$, in such a way that the hypotheses (1) and (2) of Theorem \ref{edimTheorem} are met.  Putting $w = \dim \F L_{e-t\mu}^S (\Gamma)=\dim \F L_{e-t\mu}^S(\Gamma)|_C$, we then have $$\vdim \F L_{e+s\mu}^T (\Gamma';\mu^{dst};w) \leq 0,$$ and we must only show this series is empty.  Recall that the $dst$ $\mu$-uple points lie on the curve $C$, supported on a divisor of class $dH\in \Pic C$, and that this series consists of members of $\F L_{e+s\mu}^T (\Gamma';\mu^{dst})$ that restrict to the $w$-dimensional $\F L_{e-t\mu}^S(\Gamma)|_C$.  Now to demonstrate this series is empty, we will use Theorem \ref{LimitTheorem} to iteratively specialize some of the points in $\Gamma'$ onto $C$ in such a way that $C$ appears in the base locus of the series with a total of $\max\{\mu,1\}$ times, splitting $C$ each time as it arises; these specializations preserve the virtual dimension of the series.    After splitting the copies of $C$ from the series, we have annihilated both the $dst$ $\mu$-uple points as well as the $h^0(\OO_{C}(e-t\mu))-w$ gluing conditions coming from the series on $S$.  We are thus left with a series of curves containing some fat points and some $\delta_{m,n}$-points, all aligned along $C$.  

It might appear we haven't made any progress, since instead of having a series of general fat points we still have some points in special position along $C$.  However, if there aren't too many residual schemes lying along $C$, then they will in fact be in general position on $T$ as we allow $S$ to vary.  By making sure this  is always the case, we will only have to worry about series of general fat points and $\delta_{m,n}$-points.  We can further reduce the problem to only worrying about series of fat points by using Lemma \ref{DeltaLemma}, replacing a $\delta_{m,n}$-point by either an $m$-uple point or an $(m+1)$-uple point, and checking that both of these series have the expected dimension by using either Proposition \ref{d23} or our inductive hypothesis (once $d$ is large).

\emph{Step 1: $d=4$.}  We let a quartic surface degenerate to a union $S\cup T$ of two quadric surfaces, meeting along an elliptic quartic curve $C$.  
Since $h^0(\OO_T(2))=9$ and a general collection of simple points and $\delta_{1,1}$-points impose the expected number of conditions on $H^0(\OO_T(2))$ (in fact, by Bertini's theorem, this is true for \emph{any} series on any variety), we see that a general collection of $\alpha$ fat points and $\beta$ $\delta_{m,n}$-points lying on $C$ and (in the case of the $\delta_{m,n}$-points) pointing in the same direction as $C$ will in fact be in general position on $T$ as $S$ varies if $\alpha+2\beta \leq 8$.

\emph{Case 1: $e\geq 8$.}  Take $\mu = 0$, and specialize a subcollection $\Gamma$ of our fat points onto $S$ in such a way that the expected dimension $w$ of $\F L_e^S(\Gamma)$ is between $4$ and $13$ (which can always be done since quadruple points contribute only $10$ to the degree of $\Gamma$); specialize the remaining points $\Gamma'$ onto $T$.  It follows from Proposition \ref{d23} that $\F L_e^S(\Gamma)$ has the expected dimension and that the kernel series $\F L_{e-2}^S(\Gamma)$ is empty; the only special system $\F L_{e-2}^S(\Gamma)$ for $e\geq 8$ is $\F L_6^S(4^5)$, but the expected dimension of $\F L_8^S(4^5)$ is far larger than $13$.

Now we must show that the series $\F L_e^T(\Gamma';w)$ is empty.  We must specialize points onto $C$ with a total multiplicity of at least $w$ in order to cause $C$ to appear in the base locus.  The virtual dimension of this series is nonpositive (by Theorem \ref{edimTheorem}), so the virtual dimension of $\F L_e^T(\Gamma')$ is at most $h^0(\OO_C(e))-w$.  It follows that $$\deg \Gamma' \geq w+ h^0(\OO_T(e-2))\geq 53.$$ The sum of the multiplicities of the points in $\Gamma'$ is then clearly at least $w$ (since a quadruple point has the smallest ratio of multiplicity to degree at $2/5$, the total sum of multiplicities in $\Gamma'$ is in fact at least $\frac 25 \cdot 53$, which is far larger than $w\leq 13$) so it is possible to specialize points onto $C$ one at a time to cause $C$ to split (specialize points of highest multiplicity first).  The residual schemes consist of ordinary fat points (from specializations before $C$ splits) and perhaps a single $\delta_{m,n}$-point with $m\leq 3$ (coming from applying Theorem \ref{LimitTheorem} when the final point needed to split $C$ is specialized); furthermore, there are at most $7$ such residual schemes since specializing a point of multiplicity at least $2$ contributes that amount to the total multiplicity of points along $C$, whereas specializing a simple point does not leave a residual scheme.  It follows that the residual schemes are in general position on $T$ as $S$ varies.  We are thus left with a series $\F L_{e-2}^T(\Gamma'')$ where $\Gamma''$ has some fat points and perhaps a single $\delta_{m,n}$-point.  Use Lemma $\ref{DeltaLemma}$ to consider a pair of fat point series instead.  Such a series can only be special if it equals $\F L_6^T(4^5)$.  However, this series can never arise:  if $\Gamma'$ had a $4$-uple point, then we would have specialized it onto $C$ first, leaving us with a residual $3$-uple point.

\emph{Case 2: $e=7.$} This case is similar to the preceding one.  Take $\mu = 0$, and choose $\Gamma$ so the expected dimension $w$ of $\F L_7^S(\Gamma)$ is between $5$ and $14$.  A quick check of the list of special series in Proposition \ref{d23} verifies that $\F L_7^S(\Gamma)$ is nonspecial and $\F L_5^S(\Gamma)$ is empty--the virtual dimension of $\F L_5^S(\Gamma)$ is at most $-14$, while all special series of the form $\F L_5^S(\Gamma)$ have virtual dimension at least $0$.  Specializing the remaining points $\Gamma'$ onto $T$, we estimate $$\deg \Gamma' \geq w+ h^0(\OO_T(5)) \geq 41,$$ so again the total multiplicity of points in $\Gamma'$ is at least as big as $w$.  We then specialize points in $\Gamma'$ onto $C$ to force $C$ to split, again starting with the highest multiplicity points.  The residual schemes are again in general position on $T$ as $S$ varies, with at most one $\delta_{m,n}$-point (with $m\leq 3$).  Write this series as $\F L_5^T(\Gamma'')$, where $\Gamma''$ has some fat points and perhaps a $\delta_{m,n}$-point with $m\leq 3$.  Notice that all special series $\F L_5^T(4^a,3^b,2^c)$ have $a\geq 3$ and have virtual dimension at least $0$.  Then if $\Gamma''$ has at most one quadruple point, we replace a possible $\delta_{m,n}$-point in turn by an $m$-uple point and an $(m+1)$-uple point to deduce nonspeciality by Lemma \ref{DeltaLemma}.  If $\Gamma''$ contains at least two quadruple points, then since we specialized the highest multiplicity points first it must also contain a $3$-uple point and a $\delta_{3,m}$-point with $0\leq m\leq 3$.  Replacing the $\delta_{3,m}$ point by a triple point will yield a series with two triple points, so the virtual dimension will be negative if there are three quadruple points, and thus the series is nonspecial.    On the other hand, replacing the $\delta_{3,m}$ point by a quadruple point makes the virtual dimension of the series negative, and hence also yields a nonspecial series.

\emph{Case 3: $e=6$.}  Here our strategy changes a bit.  We take $\mu = 3$, and specialize all the points onto $T$, leaving us with a series $\F L_{12}^T(\Gamma;3^{16})$ with the $16$ points where the triple points are supported forming a divisor of class $4H\in \Pic C$; the series on $S$ is the $1$-dimensional $\F L_0^S(\emptyset)$, and the kernel series on $S$ is clearly empty.

On $T$, it takes $1=h^0(\OO_C)$ condition for $C$ to appear in the base locus of $\F L_{12}^T(\Gamma;3^{16})$, a further $8=h^0(\OO_C(2))$ to force it to appear a second time, and $16= h^0(\OO_C(4))$ more to cause it to appear a third time.    Notice that $\deg \Gamma\geq 74$, so the total multiplicity of the points in $\Gamma$ is far bigger than $16+8+1=25$, and there are an abundance of points with which to split $C$ three times.   Furthermore, $\Gamma$ either
\begin{enumerate}
\item contains three quadruple points,
\item contains three triple points,
\item contains four double points, or
\item contains 11 simple points,
\end{enumerate}
since if none of the above hold then $\deg \Gamma\leq 2\cdot 10+2\cdot 6+3\cdot 3+10\cdot 1 = 51.$  Choose the \emph{first} of the four above conditions that holds, and specialize all the guaranteed points onto $C$ one at a time, splitting $C$ from the series as it appears in the base locus.  The resulting residual schemes after this specialization are encoded in Figure 1; for instance, after specializing the three quadruple points, $C$ splits twice and we are left with a $\delta_{2,2}$-point (where the first quadruple point limited onto $C$), a triple point (where the second quadruple point limited onto $C$), and a quadruple point.  Each of these four possible specializations causes $C$ to split twice.

\begin{figure}[t]
\input{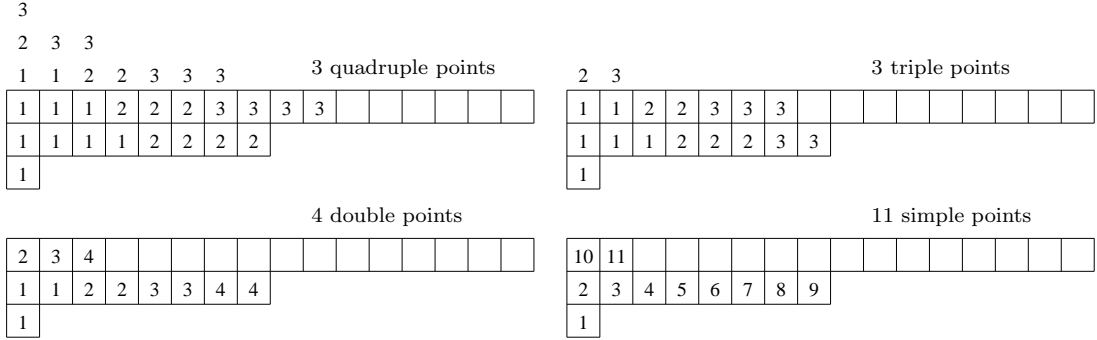}%
\caption{This diagram computes the residual schemes after specializing points onto $C$ to cause $C$ to split twice in the case $d=4$, $e=6$ of the theorem.  The boxed rows are $1$ wide, $8$ wide, and $16$ wide, corresponding to the number of conditions it takes to split $C$ one, two, and three times.  Numbers denote the order in which the corresponding points were specialized.  Numbers above the boxed rows correspond to residual schemes after $C$ splits a third time.}
\end{figure}

To force $C$ to split a third time, we can specialize points onto $C$ one at a time in an arbitrary fashion.  If we are using the first specialization, with $3$ quadruple points, then there will be at most $6$ residual schemes, at most two of which are $\delta_{m,n}$-points (one of which is the $\delta_{1,1}$-point coming from the first quadruple point that was specialized).  On the other hand, if we are using any of the other specializations, there will be at most $7$ residual schemes, at most one of which is a $\delta_{m,n}$-point. Thus in every case, the residual schemes are in general position on $T$ as $S$ varies. We then prove the nonspeciality of the residual series by replacing any $\delta_{m,n}$-points by two different fat points in turn.  The only special series $\F L_6^T(\Gamma')$ is $\F L_6^T(4^5)$, and this series clearly does not arise from this process: if our original $\Gamma$ had at least $5$ quadruple points, then we would have taken the first specialization strategy, and there would be a residual triple point.   

\emph{Case 4: $e=5$ and $e=4$.}  These cases are extremely similar to the previous case, but considerably easier.  One can take $\mu=2$ and again specialize all points onto $T$.  From there, specialize points onto $C$, starting with the points of highest multiplicity, until $C$ splits twice.  One easily checks that the residual schemes are in general position and that the resulting series is nonempty by using Lemma \ref{DeltaLemma} and Proposition \ref{d23}.  We leave the full details to the reader.

\emph{Case 5: $e=3$.}  We already showed that the series $\F L_3^4(4^2)$ is nonspecial of virtual dimension $0$ in Lemma \ref{twoQuadruple}, so we may assume there is at most one quartuple point.  We take $\mu=0$, and notice that $h^0(\OO_S(3)) = 16$.  Specialize some points $\Gamma$ onto $S$ in such a way that the virtual dimension $w$ of $\F L_3^S(\Gamma)$ is either $4$ (by specializing a combination of simple, double, and/or triple points) or $6$ (by specializing the lone quadruple point).  Then the series on $S$ is nonspecial and the kernel series is empty.  Specializing the remaining points $\Gamma'$ onto $T$, we must specialize $w$ conditions onto $C$ to force it to split.  In case $w=4$, we have $\deg \Gamma'\geq 8$, and this can clearly be done.  Likewise, if $w=6$ then $\deg \Gamma'\geq 10$, and the total multiplicity of points in $\Gamma'$ is at least $6$ (since there is no quadruple point).  In either case, the original series is nonspecial.

\emph{Case 6: $e=2$.}  We will handle this case by a more general argument in a moment.

This completes the proof for $d=4$.

\emph{Step 2: $d\geq 5$.}  First assume $e< d-1$.  In this case we degenerate a surface of degree $d$ into the union of a plane $S$ and a general surface $T$ of degree $d-1$, specializing all the fat points $\Gamma$ onto $T$; we take $\mu=0$.  The curve $C=S\cap T$ has degree $d-1$, and the hypothesis $e<d-1$ means the restriction map $H^0(\OO_S(e))\to H^0(\OO_C(e))$ is an isomorphism.  Thus $\F L_e^{T}(\Gamma)\cong \F L_e^d(\Gamma)$, and by induction on $d$ we find that $\F L_e^d(\Gamma)$ is empty unless it equals $\F L_2^d(4)$.  Note that this same argument also proves the $d=4$, $e=2$ case above, using the classification of special series on a cubic surface as a starting point for the induction.

We therefore assume that $e\geq d-1$.  In every case, we degenerate a surface of degree $d$ into a quadric $S$ and a general surface $T$ of degree $d-2$, meeting along a curve $C$ of degree $2(d-2)$.  

\emph{Case 1: $(e,d)\in \{(4,5),(5,5),(5,6)\}$.}  In each case, we take $\mu=1$ and specialize all the points $\Gamma$ onto $T$, to get a series $\F L_{e+2}^T(\Gamma;1^{2d(d-2)})$ with the $2d(d-2)$ simple points forming a divisor of class $dH$ on $C$.  The kernel series on $S$ is empty.  We have $h^0(\OO_C(e-d+2))\leq 9$, and there are plenty of points in $\Gamma$ to cause $C$ to split upon specialization. When $C$ splits we are left with a series $\F L_e^{d-2}(\Gamma'')$ with residual schemes in general position.  Replacing any $\delta_{m,n}$-points by fat points as appropriate, we conclude by induction that the series is nonspecial.

\emph{Case 2: $e\geq 6$.}  We take $\mu = 0$.  Divide our fat points into two collections $\Gamma$, $\Gamma'$, where $\Gamma$ is chosen so that the expected dimension $w$ of $\F L_e^S(\Gamma)$ is between $7$ and $16$.  By Proposition \ref{d23}, $\F L_e^S(\Gamma)$ is nonspecial.  The virtual dimension of the kernel series $\F L_{e-(d-2)}^S(\Gamma)$ is so negative that it is clearly also empty by the proposition.

Working on the surface $T$ of degree $d-2$, we must show the series $\F L_e^T(\Gamma';w)$ is empty.  The scheme $\Gamma'$ has degree at least $h^0(\OO_T(e-2))+w \geq 41$, from which it follows that the sum of the multiplicities of points in $\Gamma'$ is at least $w$ (this is the reason we had to treat the cases in Case 1 above separately).  We can thus specialize points in $\Gamma'$ onto $C$ one at a time arbitrarily to force it to split.  Since $w\leq 16$, the residual schemes will be in general position on $T$ as $S$ varies:  there are at most $8$ residual schemes, at most one of which is a $\delta_{m,n}$-point instead of a fat point.  Since $h^0(\OO_T(2))=10$, the schemes are guaranteed to be in general position.  Dealing with a potential $\delta_{m,n}$-point by using Lemma \ref{DeltaLemma}, our series is nonspecial by induction on $d$.
\end{proof}

\bibliographystyle{hsiam}

\begin{thebibliography}{10}

\bibitem{BocciMiranda}
C.~Bocci and R.~Miranda.
\newblock Topics on interpolation problems in algebraic geometry.
\newblock {\em Rend. Sem. Mat. Univ. Politec. Torino}, 62(4):279--334, 2004.

\bibitem{CC}
L.~Chiantini and  C.~Ciliberto.
\newblock Weakly defective varieties. 
\newblock {\em Trans. Amer. Math. Soc.}, 354, 151–-178, 2002. 

\bibitem{CM1}
C.~Ciliberto and R.~Miranda.
\newblock Degenerations of planar linear systems.
\newblock {\em J. Reine Angew. Math.}, 501:191--220, 1998.

\bibitem{CM2}
C.~Ciliberto and R.~Miranda.
\newblock Linear systems of plane curves with base points of equal
  multiplicity.
\newblock {\em Trans. Amer. Math. Soc.}, 352(9):4037--4050, 2000.

\bibitem{CM3} 
C.~Ciliberto and R.~Miranda.
\newblock The Segre and Harbourne-Hirschowitz conjectures, in Applications of algebraic geometry to coding theory, physics and computation.
\newblock {\em NATO Sci. Ser. II Math. Phys. Chem. 36}, Kluwer Acad. Publ., Dordrecht, 2001. 

\bibitem{Laface2}
C.~De~Volder and A.~Laface.
\newblock Degeneration of linear systems through fat points on {$K3$} surfaces.
\newblock {\em Trans. Amer. Math. Soc.}, 357(9):3673--3682 (electronic), 2005.

\bibitem{Laface}
C.~De~Volder and A.~Laface.
\newblock Linear systems on generic {$K3$} surfaces.
\newblock {\em Bull. Belg. Math. Soc. Simon Stevin}, 12(4):481--489, 2005.

\bibitem{Dumnicki}
M.~Dumnicki.
\newblock Special homogeneous linear systems on {H}irzebruch surfaces. \newblock {\em Geometriae Dedicata,} 147(1):283-311, 2010.

\bibitem{DumnickiJarnicki}
M.~Dumnicki and W.~Jarnicki.
\newblock New effective bounds on the dimension of a linear system in {$\Bbb
  P^2$}.
\newblock {\em J. Symbolic Comput.}, 42(6):621--635, 2007.

\bibitem{Evain}
L.~Evain.
\newblock Computing limit linear series with infinitesimal methods.
\newblock {\em Ann. Inst. Fourier (Grenoble)}, 57(6):1947--1974, 2007.

\bibitem{gimi2}
A.~Gimigliano.
\newblock {\em On linear systmes of plane curves}.
\newblock PhD thesis, Queen's University, Kingston, 1987.

\bibitem{Gimi}
A.~Gimigliano.
\newblock Our thin knowledge of fat points.
\newblock In {\em The {C}urves {S}eminar at {Q}ueen's, {V}ol.\ {VI}
  ({K}ingston, {ON}, 1989)}, volume~83 of {\em Queen's Papers in Pure and Appl.
  Math.}, pages Exp.\ No.\ B, 50. Queen's Univ., Kingston, ON, 1989.

\bibitem{GriffithsHarris}
P.~Griffiths and J.~Harris.
\newblock On the {N}oether-{L}efschetz theorem and some remarks on
  codimension-two cycles.
\newblock {\em Math. Ann.}, 271(1):31--51, 1985.

\bibitem{Harbourne}
B.~Harbourne.
\newblock The geometry of rational surfaces and {H}ilbert functions of points
  in the plane.
\newblock In {\em Proceedings of the 1984 {V}ancouver conference in algebraic
  geometry}, volume~6 of {\em CMS Conf. Proc.}, pages 95--111, Providence, RI,
  1986. Amer. Math. Soc.

\bibitem{Hirschowitz}
A.~Hirschowitz.
\newblock Une conjecture pour la cohomologie des diviseurs sur les surfaces
  rationnelles g\'en\'eriques.
\newblock {\em J. Reine Angew. Math.}, 397:208--213, 1989.

\bibitem{Laface3}
A.~Laface.
\newblock On linear systems of curves on rational scrolls.
\newblock Geometriae Dedicata 90 (2002), no. 1, 127--144, 2002.

\bibitem{Roe}
J.~Roe.
\newblock Limit linear systems and applications, 2006.

\bibitem{Segre}
B.~Segre.
\newblock Alcune questioni su insiemi finiti di punti in geometria algebrica.
\newblock {\em Univ. e Politec. Torino Rend. Sem. Mat.}, 20:67--85, 1960/1961.

\bibitem{Yang}
S.~Yang.
\newblock Linear systems in {${\Bbb P}^2$} with base points of bounded
  multiplicity.
\newblock {\em J. Algebraic Geom.}, 16(1):19--38, 2007.

\end{thebibliography}

\end{document}